\documentclass[reqno]{amsart}
\usepackage{amssymb,hyperref}

\newtheorem{thm}{Theorem}
\newtheorem{lem}{Lemma}
\newtheorem{cor}{Corollary}

\renewcommand{\Re}{\mathrm{Re}}

\title
[An application of Jack's lemma for the minimum point]
{An application of Jack's lemma \\
for the minimum point}

\author{Hitoshi Shiraishi}
\address{Hitoshi Shiraishi \newline
Department of Mathematics \newline
Kinki University \newline
Higashi-Osaka, Osaka 577-8502, Japan}
\email{step\_625@hotmail.com}

\date{}

\begin{document}

\begin{abstract}
For the analytic function $f(z)$,
H. Shiraishi and S. Owa (Stud. Univ. Babe\c{s}-Bolyai Math. {\bf 55}(2010), 207-211) have shown a theorem for the minimum value of $|f(z)|$.
In this paper, we discuss an application of this theorem and some corollaries.
\end{abstract}

\maketitle

\section{Introduction}

\

Let the set of $\mathbb{U}$ be the open unit disk $\{ z\in\mathbb{C} : |z|<1 \}$,
and let $\mathcal{H}[a_0,n]$ denote the class of functions $p(z)$ of the form
$$
p(z)
= a_0 +\sum^{\infty}_{k=n} a_k z^k
$$
which are analytic in $\mathbb{U}$ for some $a_0\in\mathbb{C}$ and a positive integer $n$.

\

The basic tool in the proof of our results is the following lemma due to H. Shiraishi and S. Owa \cite{d5ref0}.

\

\begin{lem} \label{d5thm1} \quad
Let $p(z)\in\mathcal{H}[a_0,n]$ with $p(z) \neq 0$ for all $z\in\mathbb{U}$.
If there exists a point $z_0\in\mathbb{U}$ such that
$$
\min_{|z| \leqq |z_0|} |p(z)|
=|p(z_0)|,
$$
then
$$
\frac{z_0 p'(z_0)}{p(z_0)}
= -m
$$
and
$$
\Re\frac{z_0 p''(z_0)}{p'(z_0)}+1
\geqq -m,
$$
where
$$
m
\geqq n \frac{|a_0-p(z_0)|^2}{|a_0|^2-|p(z_0)|^2}
\geqq n \frac{|a_0|-|p(z_0)|}{|a_0|+|p(z_0)|}.
$$
\end{lem}

\

\section{Main theorem}

\

Applying Lemma \ref{d5thm1},
we derive

\

\begin{thm} \label{p05thm1} \quad
Let the function $f(z)$ given by
$$
f(z) = a_n z^n + a_{n+l} z^{n+l} + a_{n+l+1} z^{n+l+1} + \ldots
\qquad (a_n, a_{n+l} \neq 0)
$$
be analytic in $\mathbb{U}$ and $f(z) \neq 0$ for $z\in\mathbb{U}\setminus\{ 0 \}$.
If there exists a point $z_0\in\mathbb{U}\setminus\{ 0 \}$ such that
$$
\min_{|z| \leqq |z_0|} |f(z)|
=|f(z_0)|,
$$
then
\begin{equation}
\frac{z_0 f'(z_0)}{f(z_0)}
= n-m
\label{p05thm1eq1}
\end{equation}
and
\begin{equation}
(n-m) \left( \Re\frac{z_0 f''(z_0)}{f'(z_0)}+1 \right)
\leqq (n-m)^2,
\label{p05thm1eq2}
\end{equation}
where
$$
m
\geqq l \frac{|a_n z_0^n - f(z_0)|^2}{|a_n z_0^n|^2-|f(z_0)|^2}
\geqq l \frac{|a_n z_0^n|-|f(z_0)|}{|a_n z_0^n|+|f(z_0)|}.
$$
\end{thm}

\

\begin{proof} \quad
We define the function $p(z)$ by
\begin{align*}
p(z)
& = \frac{f(z)}{z^n} \\
& = a_n + a_{l} z^{l} + a_{l+1} z^{l+1} + \ldots.
\end{align*}

Then,
$p(z)\in\mathcal{H}[a_n,l]$ and $p(0)=a_n \neq 0$.
Furthermore,
by the assumtion of the theorem,
$|p(z)|$ takes its minimum value at $z=z_0$
in the closed disk $|z| \leqq |z_0|$.
It follows from this that
$$
|p(z_0)|
= \frac{|f(z_0)|}{|z_0|^n}
= \frac{\displaystyle{\min_{|z| \leqq |z_0|}} |f(z)|}{|z_0|^n}
= \min_{|z| \leqq |z_0|} |p(z)|.
$$

Therefore,
applying Lemma \ref{d5thm1} to $p(z)$,
we observe that
$$
\frac{z_0 p'(z_0)}{p(z_0)}
= \frac{z_0 f'(z_0)}{f(z_0)} -n
= -m
$$
which shows (\ref{p05thm1eq1}) and
\begin{align*}
\Re\frac{z_0 p''(z_0)}{p'(z_0)} +1
& = \Re \left( -n-1 + \frac{\dfrac{z_0f''(z_0)}{f'(z_0)}+1-n}{1-n\dfrac{f(z_0)}{z_0f'(z_0)}} \right) +1 \\
&= -\frac{n-m}{m} \left( \Re\frac{z_0 f''(z_0)}{f'(z_0)} +1-n \right) -n \\
&\geqq -m
\end{align*}
which implies (\ref{p05thm1eq2}),
where
$$
m
\geqq l \frac{|a_n-p(z_0)|^2}{|a_n|^2-|p(z_0)|^2}
= l \frac{|a_n z_0^n - f(z_0)|^2}{|a_n z_0^n|^2-|f(z_0)|^2}
\geqq l \frac{|a_n z_0^n|-|f(z_0)|}{|a_n z_0^n|+|f(z_0)|}.
$$

This completes the assertion of Theorem \ref{p05thm1}.
\end{proof}

\

Letting $l=n$ in Theorem \ref{p05thm1},
we obtain

\

\begin{cor} \label{p05cor1} \quad
Let the function $f(z)$ given by
$$
f(z) = a_n z^n + a_{2n} z^{2n} + a_{2n+1} z^{2n+1} + \ldots
\qquad (a_n, a_{2n} \neq 0)
$$
be analytic in $\mathbb{U}$ and $f(z) \neq 0$ for $z\in\mathbb{U}\setminus\{ 0 \}$.
If there exists a point $z_0\in\mathbb{U}\setminus\{ 0 \}$ such that
$$
\min_{|z| \leqq |z_0|} |f(z)|
=|f(z_0)|,
$$
then
$$
\frac{z_0 f'(z_0)}{f(z_0)}
= n-m
\leqq 0
$$
and
$$
\Re\frac{z_0 f''(z_0)}{f'(z_0)}+1
\geqq n-m,
$$
where
$$
m
\geqq n \frac{|a_n z_0^n - f(z_0)|^2}{|a_n z_0^n|^2-|f(z_0)|^2}
\geqq n \frac{|a_n z_0^n|-|f(z_0)|}{|a_n z_0^n|+|f(z_0)|}.
$$
\end{cor}

\

Moreover,
putting $n=1$ and $a_n=1$ in Theorem \ref{p05thm1},
we get the following corollary due to M. Nunokawa and S. Owa \cite{p05ref1}.

\

\begin{cor} \label{p05cor2} \quad
Let the function $f(z)$ given by
$$
f(z) = z + a_{l+1} z^{l+1} + a_{l+2} z^{l+2} + \ldots
\qquad (a_{l+1} \neq 0)
$$
be analytic in $\mathbb{U}$ and $f(z) \neq 0$ for $z\in\mathbb{U}\setminus\{ 0 \}$.
If there exists a point $z_0\in\mathbb{U}\setminus\{ 0 \}$ such that
$$
\min_{|z| \leqq |z_0|} |f(z)|
=|f(z_0)|,
$$
then
$$
\frac{z_0 f'(z_0)}{f(z_0)}
= 1-m
\leqq 0
$$
and
$$
\Re\frac{z_0 f''(z_0)}{f'(z_0)}+1
\geqq 1-m,
$$
where
$$
m
\geqq l \frac{|z_0 - f(z_0)|^2}{|z_0|^2-|f(z_0)|^2}
\geqq l \frac{|z_0|-|f(z_0)|}{|z_0|+|f(z_0)|}.
$$
\end{cor}

\


\begin{thebibliography}{}

\bibitem{m1ref1}
I. S. Jack,
{\it Functions starlike and convex of order $\alpha$},
J. London Math. Soc. {\bf 3}(1971), 469--474.

\bibitem{m1ref2}
S. S. Miller and P. T. Mocanu,
{\it Second-order differential inequalities in the complex plane},
J. Math. Anal. Appl. {\bf 65}(1978), 289--305.

\bibitem{p05ref1}
M. Nunokawa and S. Owa,
{\it Notes on certain analytic functions},
Proc. Japan Acad. {\bf 65}(1989), 85--88.

\bibitem{d5ref0}
H. Shiraishi and S. Owa,
{\it An application of Miller and Mocanu lemma},
Stud. Univ. Babe\c{s}-Bolyai Math. {\bf 55}(2010), 207--211.

\end{thebibliography}
\end{document}